\documentstyle[12pt]{article}

\makeatletter
\renewcommand{\@biblabel}[1]{#1.}
\makeatother
\newtheorem{Sa}{\sc \hspace{0.4cm} Proposition}

\newtheorem{Th}[Sa]{\sc \hspace{0.4cm} Theorem}

\newcommand{\RR}{{\rm I\kern-0.14em R }}
\newcommand{\NN}{{\rm l\kern-0.14em N}}
\newcommand{\CC}{{\rm\raise 0.192ex\vbox{\hrule height
                  1.22ex width 0.8pt}\kern-0.29em C}}

\begin{document}
\begin{flushleft}
{\bf \Large  Pascal's Theorem and Quantum Deformation }
\end{flushleft}


\begin{flushleft}
{\large FRANK LEITENBERGER } \\
{\small \it Fachbereich Mathematik, Universit\"at Rostock,
Rostock, D-18051, Germany. \\
e-mail: frank.leitenberger@mathematik.uni-rostock.de }
\end{flushleft}


\small \normalsize

\vspace{0.2cm}

\noindent
{\small {\bf Abstract.}
By a transfer principle
Pascal's Theorem is equivalent
to a theorem about point pairs on the real line.
It appears that Pascal's Theorem is equivalent
to the vanishing of a common invariant of six quadratic forms.
Using the q-deformed invariant theory of \cite{L}
we construct corresponding quantum invariants
by a computer calculation.
}
\vspace{0.2cm}

\begin{flushleft}
\end{flushleft}

\begin{flushleft}
\end{flushleft}

\vspace{1.0cm}
\setcounter{equation}{0}

\begin{flushleft}
\bf 1. Introduction
\end{flushleft}

\noindent
Noncommutative geometry and the theory of quantum groups
present quantisations of projective spaces (cf. \cite{Z}).
Less is known
about quantum versions of fundamental relations
of points and lines in these spaces.
It is the aim of this letter to demonstrate
that the theory of Quantum groups
allows to quantise situations of elementary geometry.

We consider six points upon a conic.
They form 60 different hexagons.
For each of these hexagons the intersections
of opposite lines are in line (Pascal's Theorem).
In connection with this Theorem
we consider the following Proposition (cf. \cite{H1}).

\begin{Sa}
Let $P_a,P_b,P_c,P_d,P_e,P_f$ be six points
on the projective real line.
Let $\{A_1,B_1\}$ be a  point pair
in harmonic position to the point pairs $\{P_a,P_b\}$ and $\{P_d,P_e\}$.
Further let $\{A_2,B_2\}$ be a point pair
in harmonic position to the point pairs $\{P_b,P_c\}$ and $\{P_e,P_f\}$,
and let $\{A_3,B_3\}$ be a point pair
in harmonic position to the point pairs $\{P_c,P_d\}$ and $\{P_f,P_a\}$.
Then the three point pairs
$\{A_1,B_1\}$, $\{A_2,B_2\}$ and $\{A_3,B_3\}$ are in involution
(i.e., there is an involutive projective map,
which changes the points of each of the three point pairs).
\end{Sa}

According to the "\"Ubertragungsprinzip" of Hesse
(cf. \cite{H2})
Proposition 1 is equivalent to Pascal's Theorem:
Consider a conic $e$ with tangent $t$.
Let $P$ be a point outside $e$ and let $A,B$ be the intersections
of both tangents from $P$ on $e$ with $t$.
In this manner the points outside $e$ correspond unique to the
point pairs of $t$.
In the limit the points on $e$
correspond to (double) points of $t$.
(For simplicity we do not consider the inner points of $e$
which correspond to pairs of complex conjugated points of $t$.)

Our correspondence has the property,
that three points (outside $e$) are in line
if the corresponding point pairs on $t$
are in involution (cf. \cite{H2}).
This property allows to transfer Propositions
of two dimensional projective geometry
to Propositions of one dimensional geometry.

Furthermore let $S$ be the intersection of two lines $s_i$, $i=1,2$
intersecting $e$.
Then the 4 intersection points on $e$ correspond, respectively,
to 4 (double) points $C_i,D_i$, $i=1,2$ on $t$,
and $S$ corresponds to the point pair on $t$
in harmonic position to $\{ C_1,D_1 \}$ and $\{C_2,D_2\}$
(cf. \cite{H2}).
Now we consider six points upon $e$.
They form 60 different hexagons.
For a fixed hexagon we denote the six points
according their succession by $E_1,E_2,...,E_6$.
(The choice of the start point $E_1$
and the direction play no role.)
These points correspond, respectively,
to six (double) points $P_1,P_2, \cdots , P_6$ on $t$.
(For simplicity we restrict the consideration
to only 19 of 60 cases,
where the $E_i$ have a succession on $e$
such that the intersections $S_i$ of the lines
trough  $E_i,E_{i+1}$ and $E_{i+3},E_{i+4}$,
($i=1,2,3$; $E_7:=E_1$) are outside of $e$.
This holds for 19 of the 60 Pascal configurations).
According to this construction
Proposition 1 is equivalent to Pascal's Theorem.

In the following we give Proposition 1
an invariant theoretic form.
We describe the points $P_i$ by homogeneous coordinates $x_i,y_i$
or as zero of the linear form $(0i):=xy_i-yx_i$.
Furthermore we consider the determinants $(ij):=x_i y_j-y_i x_j$.
It is known, that two point pairs $(P_1,P_2)$ and $(P_3,P_4)$
are in harmonic position if we have for the cross ratio
\[  \frac{ \overline{P_1 P_3}\ \cdot \ \overline{P_4 P_2} }{
           \overline{P_1 P_4}\ \cdot \ \overline{P_3 P_2} }  =
    \frac{(13)(42)}{(14)(32)} = -1.  \]

Furthermore we consider point pairs $\{P_i,P_j\}$ as the zeros
of the quadratic forms $f_{ij}:=(0i)(0j)$.
Then the (unique) point pair,
which is in harmonic position to the zeros
of the two quadratic forms $f,f'$ is given by the Jacobian
\[  J_{f,f'} := \left| \begin{array}{cc}
               f_{x} & {f'}_{x}   \\
               f_{y} & {f'}_{y}
\end{array} \right|.    \]
The zeros of the three quadratic forms
$f_i=a_i x^2 + 2b_i xy + c_i y^2$, $i=1,2,3$
are in involution if we have for their combinant
\[
C_{f_1,f_2,f_3}:=
\left| \begin{array}{ccc}
               a_1 & b_1 & c_1  \\
               a_2 & b_2 & c_2  \\
               a_3 & b_3 & c_3
\end{array} \right|  = 0  \]
(cf. \cite{H2}).
Therefore Pascal's Theorem is equivalent
to a relation between quadratic forms.

\begin{Th} Let $(01),(02),...,(06)$ be six linear forms.
Let $f_1=a_1 x^2 + 2b_1 xy + c_1 y^2$ be the Jacobian
of the two quadratic forms $(01)(02)$ and $(04)(05)$.
Analogous let $f_2$ be the Jacobian  of $(02)(03)$ and $(05)(06)$,
and let $f_3$ be the Jacobian of $(03)(04)$ and $(01)(06)$.
Then we have $C_{f_1,f_2,f_3}=0$.
\end{Th}

\noindent
{\it Remark}:  The Theorem corresponds to the
associativity law  of the multiplication on the real line.

\begin{flushleft}
\bf 2. Noncommutative invariant theory
\end{flushleft}

In this section we introduce some notions from \cite{L}.
Let $q\in \CC $ with $|q|=1$ and $q\neq \pm 1$.
We replace the algebra of homogeneous coordinates $x_i,y_i$,
where $i$ is an element of an ordered index set $I$
by the noncommutative unital $*$-algebra $H_I$
given by generators $x_i,y_i$
and relations
\begin{eqnarray}
\begin{array}{lll}
  x_iy_i             &  = &      q   y_i x_i,   \\
  x_jx_i             &  = &      q^2 x_i x_j,   \\
  y_jy_i             &  = &      q^2 y_i y_j,   \\
  x_jy_i             &  = &      q   y_i x_j
                  + ( q^2 -1 )       x_i y_j,   \\
  y_jx_i             &  = &      q   x_i y_j
\end{array}
\end{eqnarray}
for $i<j$. We define the involution by
$x_i^*:=x_i$, $y_i^*:=y_i$.

$H_I$ admits the action of a Quantum group.
Let $U_q (sl(2,\RR ))$
be the unital $*$-algebra
determined by generators
$E$, $F$, $K$, $K^{-1}$
and relations
\[  K E = q E K,
 \ \ \ \ \  K F     = q^{-1} F K,
 \ \ \ \ \ [E,F]=\frac{ K^2-K^{-2} }{ q-q^{-1} },
 \ \ \ \ \  K K^{-1}= K^{-1} K   =1.   \]
We endow $U_q (sl(2,\RR ))$ with the involution
$  E^*:=F,     \ \ F^*:=E, \ \
   K^*:=K^{-1} $.
The action of $U_q(sl(2,\RR ))$ on $H_I$ is determined
if we set
\begin{eqnarray*}
\begin{array}{lll}
K 1   =  1                   ,\ \ \ \ \   &
K x_i =  q^{-\frac{1}{2}} x_i,            &
K y_i =  q^{ \frac{1}{2}} y_i,\\
E 1   =  0                   ,\ \ \ \ \   &
E x_i =  q^{\frac{1}{2}}  y_i,            &
E y_i =  0                   ,\\
F 1   =  0                   ,\ \ \ \ \   &
F x_i =  0                   ,            &
F y_i =  q^{-\frac{1}{2}} x_i
\end{array}
\end{eqnarray*}
and require that
\begin{eqnarray*}
\begin{array}{lll}
 K(ab)  &  =  &  K(a) K(b),               \\
 E(ab)  &  =  &  E(a) K(b)+K^{-1}(a)E(b), \\
 F(ab)  &  =  &  F(a) K(b)+K^{-1}(a)F(b)
\end{array}
\end{eqnarray*}
for $a,b\in H_I$. The algebra $H_I$ has the PBW-property
with respect to every order of the variables.
For further properties we refer to \cite{L}.

We say that $a \in H_I$ is an {\it invariant element}
if $Ea=Fa=0$, $Ka=a$.
The invariant elements form a subalgebra $H^{Inv}\subset H_I$.
The simplest invariant is the bracket symbol
\begin{eqnarray*}
(ij):= q^{-\frac{1}{2}} x_i y_j- q^{\frac{1}{2}} y_i x_j
\ \ \ \ \    i,j \in I.
\end{eqnarray*}
The subalgebra $H^{Inv}$
is generated by the elements $(ij)$ (cf. \cite{L}, p.94).
We have the properties
$ (ji)   = - (ij)$,
$ (ii)   =   0   $,
$ (ij)^* =   (ij)$, $\forall i,j$.
Furthermore we have for $i<j<k<l$
\[
\begin{array}{lll}
 (jk)(ij)  =  q^{ 4} (ij)(jk), \ \ \ \ \
 (ik)(ij)  =  q^{ 2} (ij)(ik), \ \ \ \ \
 (jk)(ik)  =  q^{ 2} (ik)(jk).
\end{array}
\]

In the following we assume for $H_I$ that $0\in I$.
We use the notations $x:=x_0$, $y:=y_0$.
$H_{  I \backslash \{ 0 \} }$ is the subalgebra
generated by the elements $x_i,y_i$,
$i\in I \backslash \{ 0 \}$.

We say that
$f=x^2 A+[2]_q xy B+y^2 C$ is a {\it quadratic form},
if $A,B,C\in H_{I\backslash \{ 0 \} }$ and $f\in H_I^{Inv}$.
($[2]_q$ denotes the $q$-number $q+\frac{1}{q}$.)
The second requirement corresponds to the classical requirement
for the transformation behaviour of the coefficients.
As a consequence of the PBW-Theorem
the representation of $f$
with coefficients $A,B,C$ is unique.

{\it Remark:} Alternatively, one can consider quadratic forms
$f=A_L x^2+[2]B_L xy+C_L y^2$ with left coefficients.
If $0$ is the minimal element of $I$,
then one can show that the right and left coefficients
differ only by certain powers of $q$.

In the following we consider the  quadratic forms
\[  f_{ij} := q(0i)(0j) = \frac{1}{q} (0j)(0i)  \ \ \ \ \
{\rm with} \ \ \ \ \   i<j.   \]
We obtain
\begin{eqnarray}
A =     A_L = q y_1 y_2,\ \ \ \ \
B = q^2 B_L = \frac{-q^2 x_1 y_2-q  y_1 x_2}{q+q^{-1}},\ \ \ \ \
C = q^4 C_L =   x_1 x_2.
\end{eqnarray}

Now we consider a set of quadratic forms
$f',f'',...$ with coefficients ${A'}_i,{A''}_i,...$.
We say that a polynomial expression
\[ I_{f',f'',...} = I_{f',f'',...}({A'}_i,{A''}_i,...)   \]
is a {\it common invariant}
if $I_{f',f'',...}  \in H_{I \backslash \{ 0 \} }^{Inv}$.
We say that a polynomial expression
\[ K_{f',f'',...} = K_{f',f'',...}(x,y,{A'}_i,{A''}_i,...) \]
is a {\it common quadratic covariant}
if $K_{f',f'',...}$ is a quadratic form.

In \cite{L} we have constructed invariants and covariants
by the symbolic method of Clebsch and Gordan.
For two quadratic forms $f, f'$
the {\it Jacobian}  is given by the symbol $(01)(02)(12)$.
We obtain by the symbolic method
\[ J_{f,f'} =
x^2 {  K} + xy {  L} + y^2  {  M}  \]
with
\begin{eqnarray}
\begin{array}{lll}
 {  K} = q^{6} K_L & = & - q^7 {  B A}' + q^9 {  A B}', \\
 {  L} = q^{8} L_L & = & - q^7 {  C A}' - q^6 {  B B}'
        + q^{10} {  B B}' + q^7 {  A C}',           \\
 {  M} = q^{10}M_L & = & - q^7 {  C B}' + q^9 {  B C}'.
\end{array}
\end{eqnarray}

The {\it combinant} of three quadratic forms $f,f',f''$
is given by the symbol $-\frac{1}{q^4} (12)(13)(23)$. We have
\[  C_{f,f',f''} =
 +   q^\frac{  -5}{2}{  A} {  B}' {  C}''
 -   q^\frac{  -9}{2}{  A} {  C}' {  B}''
 -   q^\frac{  -9}{2}{  B} {  A}' {  C}''
 +   q^\frac{  -9}{2}{  B} {  C}' {  A}''
 +   q^\frac{  -9}{2}{  C} {  A}' {  B}''
 -   q^\frac{ -13}{2}{  C} {  B}' {  A}'' \]
\[-(q^\frac{-7}{2}-q^\frac{-15}{2}){ B}{ B}'{ B}''.\]

In general for $q\neq 1$ the Jacobian
and the combinant do not satisfy
the classical  relations
\[
\begin{array}{c}
  J_{f,f'}     = - J_{f,f'},      \\
  C_{f,f',f''} = - C_{f,f'',f'} = - C_{f',f,f''} =
  C_{f',f'',f} =   C_{f'',f,f'} = - C_{f'',f',f}.
\end{array}
\]

\begin{flushleft}
\bf 3. Pascal's Theorem and q-deformation.
\end{flushleft}

In the following we suppose
that the linear order of the index set $I$
and the succession of the points
on the hexagon coincide (cf. the Remark below).
We have seen above,
that in the classical case the relation of invariants
\begin{eqnarray}
C_{J_{f_{12},f_{45}},J_{f_{23},f_{56}},J_{f_{34},f_{16}}} = 0
\end{eqnarray}
is equivalent to Pascal's Theorem.

One can change the arguments of $C$ and $J$ in the expression
$C_{J_{f_{12},f_{45}},J_{f_{23},f_{56}},J_{f_{34},f_{16}}}$
such that one obtains 48 expressions.
For $q\neq 1$
in general these 48 expressions do not vanish,
they are different
and they do not coincide after a possible change of the sign.
Furthermore because of (2),(3) there arise no new invariants
using left coefficients.

Therefore for a quantum analog of Theorem 2
we look for a linear combination of these 48 expressions.
In view of the relations
\[  J_{f_{45},f_{12}} = - q^4 J_{f_{12},f_{45}},  \]
\[  J_{f_{56},f_{23}} = - q^4 J_{f_{23},f_{56}}  \]
it is sufficient to consider
a linear combination of only 12 expressions instead of 48.
Therefore for an analog of equation (4) we consider
the relation

\begin{eqnarray*}
C_{a_1,...,a_{12}} & := & \ \
 a_1 C_{J_{f_{12},f_{45}},J_{f_{23},f_{56}},J_{f_{34},f_{16}}} +
 a_2 C_{J_{f_{23},f_{56}},J_{f_{34},f_{16}},J_{f_{12},f_{45}}} +
 a_3 C_{J_{f_{34},f_{16}},J_{f_{12},f_{45}},J_{f_{23},f_{56}}} \\
& & -a_4 C_{J_{f_{12},f_{45}},J_{f_{34},f_{16}},J_{f_{23},f_{56}}} -
 a_5 C_{J_{f_{34},f_{16}},J_{f_{23},f_{56}},J_{f_{12},f_{45}}} -
 a_6 C_{J_{f_{23},f_{56}},J_{f_{12},f_{45}},J_{f_{34},f_{16}}} \\
& & -a_7 C_{J_{f_{12},f_{45}},J_{f_{23},f_{56}},J_{f_{16},f_{34}}} -
 a_8C_{J_{f_{23},f_{56}},J_{f_{16},f_{34}},J_{f_{12},f_{45}}} -
 a_9C_{J_{f_{16},f_{34}},J_{f_{12},f_{45}},J_{f_{23},f_{56}}} \\
& & +a_{10}C_{J_{f_{12},f_{45}},J_{f_{16},f_{34}},J_{f_{23},f_{56}}} +
 a_{11}C_{J_{f_{16},f_{34}},J_{f_{23},f_{56}},J_{f_{12},f_{45}}} +
 a_{12}C_{J_{f_{23},f_{56}},J_{f_{12},f_{45}},J_{f_{16},f_{34}}} = 0.
\end{eqnarray*}
We obtain the following result.
\begin{Th}
We have the relation of invariants
\[  C_{a_1,...,a_{12}}=0 \]
if and only if the $a_i$ are given by
\begin{eqnarray}
\begin{array}{lll}
 a_1=a_9    & = & \ \ \alpha          \\
 a_2=a_8    & = & \ \ \beta           \\
 a_3=a_7    & = & -q^2 \frac{2q^8-5q^6+2q^4+2q^2-2}{
                  q^{12}-4q^{10}+4q^8-6q^4+6q^2-2}\alpha   \\
 a_4=a_{10} & = & \ \ 2q^2 \frac{q^{16}-3q^{14}+
                  2q^{12}+3q^{10}-6q^{8}+3q^6+2q^4-3q^2+1}{
                  q^{12}-4q^{10}+4q^8-6q^4+6q^2-2}\alpha + q^{14} \beta  \\
 a_5=a_{12} & = & -q^{-4}\frac{2q^8-2q^6-2q^4+5q^2-2}{
                  q^{12}-4q^{10}+4q^{8}-6q^4+6q^2-2}\alpha      \\
 a_6=a_{11} & = & \ \ q^{-6}\frac{2q^{12}-6q^{10}+6q^8-4q^4+4q^2-1}{
                  q^{12}-4q^{10}+4q^{8}-6q^{4}+6q^2-2}\alpha
\end{array}
\end{eqnarray}
with $\alpha, \beta \in \CC $.
\end{Th}
We set $C_{\alpha,\beta }:=C_{a_1,...,a_{12}}$.
Then we have $C_{\alpha,\beta }= \alpha C_{1,0} + \beta C_{0,1}$.

\noindent
{\it Proof:}
Applying (1) as reduction rules
we can give $C_{a_1,...,a_{12}}$ the form
\begin{eqnarray}
\sum_{\tiny \begin{array}{{c}}
0\leq i_1,...,i_6 \leq 2 \\ i_1+i_2+...+i_6=6 \end{array} }
a_{i_1,...,i_6} x_1^{i_1} y_1^{2-i_1} ... x_6^{i_6}  y_6^{2-i_6}
\end{eqnarray}
by a computer computation.

The coefficients  $a_{i_1,...,i_6}$ are linear in the $a_i$.
Therefore the equation $C_{a_1,...,a_{12}}=0$ is equivalent
to the linear equation system $a_{i_1,...,i_6}=0$
with 141 equations for 12 unknowns.
The equation system was solved by a computer program.
We found that the rank is 10
and the above two parameter solution (5).$\bullet$

\noindent
{\it Remark:}
We consider the classical case $q=1$ for $C_{1,0}$ and $C_{0,1}$.
The relation
\begin{eqnarray*}
C_{0,1 } = &
C_{J_{f_{23},f_{56}},J_{f_{34},f_{16}},J_{f_{12},f_{45}}} -
q^{14}C_{J_{f_{12},f_{45}},J_{f_{34},f_{16}},J_{f_{23},f_{56}}}
& \\ & -
C_{J_{f_{23},f_{56}},J_{f_{16},f_{34}},J_{f_{12},f_{45}}} +
q^{14}C_{J_{f_{12},f_{45}},J_{f_{16},f_{34}},J_{f_{23},f_{56}}} &
= 0   \end{eqnarray*}
reduces for $q=1$ to
$4\ C_{J_{f_{23},f_{56}},J_{f_{34},f_{16}},J_{f_{12},f_{45}}}=0$
and we obtain again Theorem 2.
But for $q=1$ the relation
\[  C_{1 ,0 } =  0 \]
reduces to the trivial relation $0=0$
because of
\[  a_1=a_5=a_9=a_{12}= 1, \ \ \ \ \ \
    a_3=a_6=a_7=a_{11}=-1, \ \ \ \ \ \
    a_2=a_4=a_8=a_{10}= 0  \]
for $\alpha =0,\beta=1, q=1$.

\noindent {\it Remark:}
We have supposed that the order of $I$ corresponds
to the succesion of the six points on the hexagon.
This is only one of 6! possibilities.
It appears that we obtain
the same above mentioned 48 invariants
in two cases if both arise by cyclic permutations
or by inversion of the direction.
Therefore it is sufficient to consider
$\frac{6!}{2 \cdot 6} = 60$ different cases.
For example for the hexagon $(125643)$
one obtains a three parameter solution $C_{\alpha,\beta,\gamma }$
of (6).
On the other side there is no nontrivial solution
of (6) for the hexagon $(123465)$.

\end{document}